\documentclass[11pt]{article}
\usepackage{amscd}
\usepackage{verbatim}
\usepackage{amssymb}
\usepackage{amsmath}
\usepackage{amsthm}
\usepackage{graphics,epsfig}

\begin{document}
\vskip 0 true cm \flushbottom

\begin{center}
\vspace{24pt} { \large \bf The Ricci flow of asymptotically
hyperbolic mass and applications} \\
\vspace{30pt}
{\bf T Balehowsky}\footnote{balehows@ualberta.ca} and
{\bf E Woolgar}\footnote{ewoolgar@math.ualberta.ca} 

\vspace{24pt} 
{\footnotesize Dept of Mathematical and Statistical
Sciences,
University of Alberta,\\
Edmonton, AB, Canada T6G 2G1.}
\end{center}
\date{\today}
\bigskip

\begin{center}
{\bf Abstract}
\end{center}

\noindent We consider the evolution of the asymptotically hyperbolic
mass under the curvature-normalized Ricci flow of asymptotically
hyperbolic, conformally compactifiable manifolds. In contrast to
asymptotically flat manifolds, for which ADM mass is constant during
Ricci flow, we show that the mass of an asymptotically hyperbolic
manifold of dimension $n\ge 3$ decays smoothly to zero exponentially
in the flow time. From this, we obtain a no-breathers theorem and a
Ricci flow based, modified proof of the scalar curvature rigidity of
zero-mass asymptotically hyperbolic manifolds. We argue that the
nonconstant time evolution of the asymptotically hyperbolic mass is
natural in light of a conjecture of Horowitz and Myers, and is a
test of that conjecture. Finally, we use a simple parabolic scaling
argument to produce a heuristic ``derivation'' of the constancy of
ADM mass under asymptotically flat Ricci flow, starting from our
decay formula for the asymptotically hyperbolic mass under the
curvature-normalized flow.


\setcounter{equation}{0}
\newpage

\section{Introduction}
\setcounter{equation}{0}

\noindent The positive energy theorem, also called the positive mass
theorem, asserts that if a complete asymptotically flat manifold has
nonnegative scalar curvature, it has nonnegative
Arnowitt-Deser-Misner (ADM) mass. One of the most interesting
aspects of this theorem is the variety of methods by which it has
been proved, each method rich in mathematical content or physical
insight and often in both. There are, for example, the Schoen-Yau
minimal surface proof \cite{SY}, Witten's spinorial method
\cite{Witten}, the spacetime fastest curves approach \cite{PSW},
Geroch's inverse mean curvature flow \cite{Geroch, HI, Bray}, and
even a method based on 3-manifold geometrization \cite{Reiris}.
(Most methods of proof also require a restriction on the dimension
$n$. The Witten method has no dimension restriction, but requires
that the manifold be spin. The work of \cite{Lohkamp} has no
dimensional or topological restriction.)

A second statement, usually referred to as the {\it rigidity}, says
that if in addition the mass is zero, then the metric is the flat
metric on ${\mathbb R}^n$. What is usually shown is that zero mass
implies Ricci flatness. While this suffices if $n=3$, for $n>3$ one
can invoke the splitting theorem to show that the unique Ricci flat
and asymptotically flat complete metric is the flat Euclidean metric
on ${\mathbb R}^n$ . Both the positivity theorem and the rigidity
statement have counterparts for asymptotically anti-de Sitter
spacetimes \cite{AbbottDeser} and for asymptotically hyperbolic
Riemannian manifolds \cite{MinOo, AnderssonDahl, Wang, CH, ACG}.

Recently, Haslhofer \cite{Haslhofer} has shown that rigidity can be
proved on asymptotically flat manifolds by using the properties of
Ricci flow. His argument depends on two important properties of
Ricci flow. These are that, during the flow of an asymptotically
flat metric, (i) the ADM mass remains constant \cite{DM, OW} and
(ii) the scalar curvature increases monotonically. Both these facts
are easy to show, once it is known that the flow of an
asymptotically flat metric remains asymptotically flat. Then the
Ricci flow beginning at a zero mass, zero scalar curvature metric
that is {\it not} Ricci flat, and thus not a fixed point of the
flow, will produce a zero-mass metric with scalar curvature
everywhere nonnegative and somewhere positive. Once this is done a
conformal transformation can be found to return the scalar curvature
to zero, strictly lowering the mass and therefore rendering it
negative, while still maintaining asymptotic flatness. However, this
would violate the positive energy theorem. Therefore, the initial
Ricci curvature had to be zero.

Very recently, Bahuaud \cite{Bahuaud} has shown that the Ricci flow
of conformally compactifiable
metrics always exists for some time interval during which the
evolving metric remains conformally compactifiable. Qing, Shi, and
Wu \cite{QSW} have studied the long-time existence and convergence
of conformally compactifiable, asymptotically hyperbolic metrics.
The flow in this case is so-called normalized flow, given by
\begin{equation}
\frac{\partial g_{ij}}{\partial t} = -2 \left (
R_{ij}+(n-1)g_{ij}\right ) =:-2E_{ij}\ . \label{eq1.1}
\end{equation}

Under various appropriate asymptotic conditions, noncompact manifolds
admit a definition of mass. Intrinsic geometric flows can sometimes
preserve these asymptotic conditions and, when they do, a natural
problem is to determine the evolution of mass during the flow; see
for example \cite{DM, OW, CZ, LQZ}.

Though the ADM mass in the asymptotically flat case is constant
under Ricci flow \cite{DM, OW}, there are good reasons to think that
mass in the asymptotically hyperbolic case will not always remain
constant during the flow (\ref{eq1.1}). One such reason is that, if
it were to remain always constant, then a version of the rigidity
argument above could be used to falsify an otherwise quite plausible
form of positive energy conjecture of Horowitz and Myers \cite{HM},
a conjecture which has consequences for physics, which we now
briefly describe.

For $\dim M\ge 3$, consider the family of metrics on $M$ given by
\begin{equation}
ds^2=\frac{dr^2}{r^2\left ( 1-\frac{1}{r^n}\right )}
+r^2\left [ \left ( 1-\frac{1}{r^n}\right ) d\xi^2
+\sum_{i=3}^{n} d\theta_i^2\right ] \ , \label{eq1.2}
\end{equation}
with $r\in [1,\infty)$, $\xi\in [0,4\pi/n]$, and $\theta_i\in
[0,a_i]$ where $0<a_3\le\dots\le a_n$. Then (\ref{eq1.2}) yields a
family (parametrized by the $a_i$) of smooth metrics on ${\mathbb
R}^2\times T^{n-2}$, where $T^{n-2}$ is an $(n-2)$-torus. (A parameter, sometimes denoted $r_0$ or $M$, often appears in
descriptions of the metric (\ref{eq1.2}) \cite{HM}, but has no significance and can be removed by rescaling the coordinates. Another parameter, $\ell$, the {\it radius of curvature at infinity}, also sometimes appears but can be removed by homothetic rescaling.) These
metrics are asymptotically (locally) hyperbolic and have scalar curvature
$R=-n(n-1)$. We can obtain an Einstein metric from (\ref{eq1.2}) by
adding an extra dimension, say with coordinate $\tau$, and adding
either $-r^2d\tau^2$ with $\tau\in {\mathbb R}$ to (\ref{eq1.2}) to
obtain a Lorentzian metric which is sometimes called the {\it AdS
soliton}, or by adding $r^2d\tau^2$ with $\tau\in S^1$ to
(\ref{eq1.2}) to obtain a Riemannian metric that has been called a
{\it toric black hole} \cite{Anderson}. To avoid confusion with
Ricci (and other) solitons and to distinguish (\ref{eq1.2}) from the
AdS soliton Einstein metric from which (\ref{eq1.2}) is induced on a
slice, we will refer to the metric (\ref{eq1.2}) by the term {\it
Horowitz-Myers geon}. (In physics, {\it geon} loosely connotes a nonsingular, stable,
localized concentration of curvature sometimes, but not
necessarily, associated with nontrivial topology. The metric
(\ref{eq1.2}) is {\it not} a Ricci soliton.)

A strong formulation of the Horowitz-Myers conjecture \cite{HM} is
that each member of this family minimizes the hyperbolic mass
amongst all metrics that asymptote to it at large $r$ and have
scalar curvature $R \ge -n(n-1)$. Interestingly, the mass of
(\ref{eq1.2}) is {\em negative}: it is
$-\frac{4\pi}{n}\prod_{i=3}^na_i$ (or
$-\frac{1}{4nG}\prod_{i=3}^na_i$ in the physics normalization
achieved by multiplying by $\frac{1}{16\pi G}$ where $G$ is Newton's
constant). Hence, the Horowitz-Myers conjecture plays the role of a
``positive'' mass conjecture for this class of asymptotic
structures, if ``positive'' is interpreted to mean $\ge
-\frac{1}{4nG}\prod_{i=3}^na_i$. The conjecture arises from the
hope, as yet unrealized, that the AdS/CFT correspondence in physics
could be extended to field theories and geometries that do not have
supersymmetry. Accepted ``low-energy'' physical theories (e.g., QCD,
the theory of the strong nuclear force) are not supersymmetric, and
if AdS/CFT could be extended to them, it would provide a new tool to
study these theories nonperturbatively (though QCD with massive
quarks is also not conformal). Horowitz and Myers argue that the
veracity of their conjecture would follow if a nonsupersymmetric
AdS/CFT correspondence were to exist. Hence, the conjecture provides
a nice test of these ideas.

The Horowitz-Myers conjecture can be used to predict that the mass
strictly increases, at least initially, along the solution of
(\ref{eq1.1}) that develops from initial data (\ref{eq1.2}). To see
how, note that the metric (\ref{eq1.2}) is not Einstein, so it is
not a fixed point of (\ref{eq1.1}). The scalar curvature of
(\ref{eq1.2}) is $-n(n-1)$ everywhere, and it is easy to show that
under the flow it will remain $\ge -n(n-1)$ and will become
$>-n(n-1)$ somewhere. We can stop the flow at some time $t>0$ and
return the scalar curvature to $-n(n-1)$ by a conformal
transformation found by solving the Yamabe equation, which can
always be done in this case \cite{AM}. This conformal transformation
lowers the mass. If the mass had remained constant (or had
decreased) under the flow, then the combination of flow followed by
a conformal transformation would produce a metric of mass less than
the soliton's mass, violating the conjecture of Horowitz and Myers.

We will in fact confirm this prediction, by showing:

\bigskip
\noindent{\bf Theorem 1.1.} {\sl Let $g(t)$, $t\in [0,T)$, be an
asymptotically hyperbolic solution of (\ref{eq1.1}) developing from
a metric $g(0)=g_0$ of mass $m_0$ on a manifold $M$ with $\dim M\ge
2$. Then the mass $m(t)$ of $g(t)$ obeys}
\begin{equation}
m(t)=m_0e^{-(n-2)t}\ . \label{eq1.3}
\end{equation}
\bigskip

\noindent For $g_0$ a Horowitz-Myers geon, then
$m_0=-\frac{4\pi}{n}\prod_{i=3}^na_i<0$ and $n\ge 3$, so from
(\ref{eq1.3}) we see that $m(t)$ is strictly monotonic increasing,
as predicted from the Horowitz-Myers conjecture. This also helps to
explain why the soliton metrics have negative, rather than zero,
mass, for if the mass were zero initially, it would remain so under
the flow. Then a conformal transformation could be found to perturb
the metric to one of negative mass, preserving $R+n(n+1)\ge 0$. A
complete, zero mass, $R+n(n+1)=0$ metric that is not Ricci flat
cannot minimize mass.

Though formula (\ref{eq1.3}) may stand in contrast to the constancy
of ADM mass in asymptotically flat Ricci flow, these two situations
are actually connected by a heuristic argument based on parabolic
scaling. We outline this argument in Section 6 and use it to provide
a heuristic derivation, starting from (\ref{eq1.3}), of the
constancy of ADM mass during asymptotically flat Ricci flow.

Various results, including some  already known and some
generalizations thereof, follow from Theorem 1.1. One that follows
immediately is

\bigskip
\noindent{\bf Corollary 1.2.} {\sl Let $g(t)$, $t\in [0,T)$, be an
asymptotically hyperbolic solution of (\ref{eq1.1}) developing from
a metric $g(0)=g_0$ of mass $m_0$ on a manifold $M$ with $\dim M\ge
3$. Then there exist times $0\le t_1<t_2< T$ such that
$m(t_1)=m(t_2)$ iff $m(t)=m_0=0$ for all $t\in [0,T)$.}
\bigskip

So-called {\it breather solutions} of $(\ref{eq1.1})$ are solutions
that are periodic up to a diffeomorphism. That is, a flow is a
breather if there are times $t_1\neq t_2$ during the flow and a
diffeomorphism $\varphi$ such that $g(t_2)=\varphi^* g(t_1)$. A
breather of the form $g(t)=\varphi^*_tg(0)$ for all $t$ is called a
{\it soliton} (solitons of (\ref{eq1.1}) are generally referred to
as {\it expanding Ricci solitons}). Einstein metrics may be regarded
as solitons of (\ref{eq1.1}) for which $\varphi_t$ is independent of
$t$.

\bigskip
\noindent{\bf Corollary 1.3 (No Massive Breathers/Solitons).} {\sl
Let $g(t)$ be as in Corollary 1.2. Let $0\le t_1<t_2<T$ be such that
$g(t_2)=\varphi^*_{t_2,t_1}g(t_1)$, where $\varphi_{t_2,t_1}-{\rm
id}\in {\cal O}(x^{\tau})$, $\tau >\frac12 \dim M$, and $x$ is a defining
function (see section 2).
Then $m(t)=m_0=0$ for all $t\in [0,T)$.}
\bigskip

\noindent{\bf Proof.} If the diffeomorphism $\varphi_{t_2,t_1}$
obeys $\varphi_{t_2,t_1}-{\rm id}\in {\cal O}(x^{\tau})$, then
$m(\varphi^*_{t_2,t_1}g(t_1))=m(g(t_1))$ by Theorem 3.4 of
\cite{Herzlich} (or Theorem 2.3 of \cite{CN};
see also Theorem 2.3 of \cite{CH}). For a breather, we
have  $m(\varphi^*_{t_2,t_1}g(t_1))=m(g(t_2))$ and thus
$m(g(t_2))=m(g(t_1))$ for some $t_1\neq t_2$. Then $m_0=0=m(t)$ by
Corollary 1.2. \qed
\bigskip

\noindent It follows that asymptotically hyperbolic Einstein
manifolds cannot have nonzero mass, as was first proved in section 5
of \cite{AnderssonDahl}.

Our results do not preclude Einstein metrics, solitons and, more
generally, breathers of undefined mass.

Asymptotically hyperbolic manifolds can be assigned a {\it
boundary-at-infinity} $\partial_{\infty}M$; see Definition 2.1
below. For example, the Horowitz-Myers geon has
$\partial_{\infty}M\simeq T^{n-1}$, while standard hyperbolic space
has $\partial_{\infty}M\simeq S^{n-1}$. For
$\partial_{\infty}M\simeq S^{n-1}$, we have \cite{MinOo,
AnderssonDahl, Wang, CH, ACG}

\bigskip
\noindent{\bf Proposition 1.4 (Rigidity).} {\sl Let $3\le n=\dim M
\le 6$. Let $M$ admit a class ${\cal G}$ of metrics whose elements
are asymptotically hyperbolic with $\partial_{\infty}M\simeq
S^{n-1}$ and $E[g]:=R+n(n-1)\ge 0$. If $M$ is not spin, then further
restrict ${\cal G}$ to those metrics whose mass aspect function (see
Definition 2.3) is of semi-definite sign. If $m[g]=0$ for some
$g\in{\cal G}$ then $(M,g)$ is isometric to standard hyperbolic
space.}
\bigskip

In fact, as with the asymptotically flat case \cite{Haslhofer}, the
rigidity theorem can be shown to be a consequence of the behaviour
of mass under the flow; that is, it is a corollary of Theorem 1.1.
This is shown in section 5, using a modified version of the proof
given in \cite{ACG}. We replace the variation of the metric used
therein by one based on the flow (\ref{eq1.1}).

Finally, from this we obtain

\bigskip
\noindent{\bf Corollary 1.5.} {\sl Let $3\le n=\dim M \le 6$, and
let $g(t_2)=\varphi^*_{t_2,t_1}g(t_1)\in{\cal G}$, with ${\cal G}$
as in Proposition 1.4. and $g$ and $\varphi_{t_2,t_1}$ as in
Corollary 1.3. Then $g$ is isometric to standard hyperbolic space.}
\bigskip

\noindent As a special case, there are no nontrivial steady solitons
of the flow (\ref{eq1.1}) with $R\ge -n(n+1)$ and
boundary-at-infinity $\partial_{\infty}M\simeq S^{n-1}$. This has
already been shown in \cite{FLW} by a different method, in the case
of $\partial_{\infty}M$ connected but of otherwise arbitrary
topology. This result is applicable to a technique to find Einstein
metrics using numerical Ricci flow, developed principally by Wiseman
(see \cite{Wiseman} and citations therein) and collaborators.

This paper is organized as follows. In section 2, we describe
asymptotically hyperbolic manifolds and mass formulae. In section 3,
we consider the normalized Ricci flow of an asymptotically
hyperbolic manifold. In section 4, we prove Theorem 1.1. In section
5, we prove Proposition 1.4. In section 6, we provide a heuristic
argument that, starting from (\ref{eq1.3}), one can predict that the
mass of an asymptotically flat metric will remain constant during
Ricci flow. An appendix discusses the equivalence of our version of
Wang's mass formula and the Chru\'sciel-Herzlich mass formula. A
second appendix outlines an alternative derivation, where we compute
the behaviour of the mass under Ricci-DeTurck flow and then pass to
the Ricci flow result by a pullback.

Our index convention for the curvature tensor is such that
$R^i{}_{jkl}x^ky^lz^j := \big ( \nabla_x\nabla_y z -
\nabla_y\nabla_x z -\nabla_{[x,y]}z\big )^i$. We maintain the
positions of indices on the Riemann tensor when raising or lowering,
so for example $R_{ijkl}:=g_{im}R^m{}_{jkl}$. We define the
Laplacian $\Delta T$ on a tensor $T$ with components $T_{i\dots
j}^{k\dots l}$ in a coordinate basis so that the components of
$\Delta T$ are $\Delta T_{i\dots j}^{k\dots l}=g^{pq}\nabla_p \left
( \nabla_q T_{i\dots j}^{k\dots l}\right )$.

\bigskip
\noindent{\bf Acknowledgements.} The authors are grateful to Eric
Bahuaud for providing helpful comments based on careful readings of
two drafts, and to him and Rafe Mazzeo for discussions. We thank an
anonymous referee for many helpful suggestions and corrections. This
work was supported by an NSERC Discovery Grant to EW.

\section{Asymptotically hyperbolic manifolds}
\setcounter{equation}{0}

\noindent There are two approaches to the question of what it means
for a manifold to be asymptotic to a hyperbolic manifold; see for
example \cite{Herzlich}. One is based on the Penrose conformal
approach first employed in general relativity, while the other uses
charts and fall-off conditions. The two approaches are commensurate,
in that conformal compactification implies that curvature and its
derivatives obey some mild decay conditions.

\bigskip
\noindent {\bf Definition 2.1.} A complete manifold $(M,g)$ is
called (smoothly) {\it conformally compactifiable} if there is a
manifold-with-boundary $({\tilde M},{\tilde g})$, a function
$\rho\in C^{\infty}({\tilde M})$ obeying $\rho(p)=0 \Leftrightarrow
p\in
\partial {\tilde M}$ and $d\rho\vert_p\neq 0$ whenever
$p\in \partial {\tilde M}$, and a smooth diffeomorphism $\psi$ from
the interior of ${\bar M}$ onto $M$ such that $\rho^2\psi^*g={\tilde
g}$ is a smooth metric on ${\tilde M}$. The function $\rho$ is
called the {\it defining function} for the {\it
boundary-at-infinity} $\partial {\tilde M}=: \partial_{\infty} M$.
\bigskip

\noindent{\bf Definition 2.2.} If $(M,g)$ is conformally
compactifiable with defining function $\rho$ such that
\begin{itemize}
\item the metric on $\partial_{\infty} M$ induced by ${\tilde
    g}$ has constant sectional curvature $k$, where either $k=0$
    or $k=1$ (or $k=-1$, but we will not treat that case in the
    sequel), and
\item $\rho$ obeys
\begin{equation}
    {\widetilde {\vert d\rho \vert}}^2:={\tilde
g}^{ij}\partial_i\rho\partial_j\rho=1+{\cal O}(\rho)\ .
    \label{eq2.1}
\end{equation}
\end{itemize}
then we call $(M,g)$ {\it asymptotically hyperbolic}.

\bigskip
\noindent Our definition is only local in that $\partial_{\infty}M$ can have
the topology of a spherical space form (for $k=1$) or a flat
manifold (for $k=0$). A standard calculation shows that (\ref{eq2.1})
implies that
\begin{equation}
\left \vert E[g] \right \vert \equiv \left \vert {\rm Ric}[g]
+(n-1)g\right \vert ={\cal O}(\rho)  \label{eq2.2}
\end{equation}
whenever ${\tilde g}\in C^2({\tilde M})$ (here
$|T_{ij}|^2:=g^{ik}g^{jl}T_{ij}T_{kl}$).

We will consider a much more limited class of metrics, admitting a
definition of mass. We first consider the $k=1$ case, with
$\partial_{\infty}M\cong S^{n-1}$ with the canonical (round) metric
${\check g}:= g(S^{n-1},{\rm can})$ (one can also consider spherical
space forms). This case was treated by Wang \cite{Wang}.
Specifically, Wang treats metrics which have an expansion
\begin{equation} g={\rm
csch}^2(x)\left [ dx^2+{\check g} +x^{n}\kappa/n +{\cal
O}(x^{n+1})\right ] \ , \label{eq2.3}
\end{equation}
where $x:={\rm arcsinh}(\rho)$, $\kappa$ is a tensor on
$\partial_{\infty}M$ (sometimes called the {\it mass aspect
tensor}). Then Wang defines the {\it mass} to be
\begin{equation}
m\equiv m[g]:=\int_{S^{n-1}}{\check g}^{AB}\kappa_{AB}d{\check \mu}
\ , \label{eq2.4}
\end{equation}
with $d{\check \mu}$ the volume element defined by ${\check g}$.
This definition agrees with a definition given by Chru\'sciel and
Herzlich \cite{CH} (see Appendix A). However, the
Chru\'sciel-Herzlich definition has the virtue of being formulated
in a general way that includes the cases where $\partial_{\infty} M$
is a spherical space form or a closed flat ($k=0$) or hyperbolic
($k=-1$) manifold. We will need to treat the $k=0$ case, so we will
extend the Wang definition to include it.

\bigskip
\noindent{\bf Definition 2.3.} Consider the set of asymptotically
hyperbolic metrics for which, in a collar neighbourhood of
$\partial_{\infty}M$ coordinatized by $(x^i)=(x,y^A)$,
the metric can be written as
\begin{eqnarray}
g&=&\frac{1}{\rho_{(k)}^2(x)}\left [ dx^2+g_{(k)} +x^n\kappa/n
+{\cal O}(x^{n+1})\right ] \ , \label{eq2.5}\\
\rho_{(k)}(x)&=&
\begin{cases}{\rm sinh}(x)\ ,&{\rm if\ }k=1,\nonumber\\
x\ ,&{\rm if\ }k=0.\nonumber
\end{cases}
\end{eqnarray}
Here $g_{(k)}$ denotes a metric on $\partial_{\infty}M$ of constant
sectional curvature $k=0$ or $k=1$, $\kappa=\kappa_{AB}dy^Ady^B$ is
a symmetric tensor on level sets of $x$, and ${\cal O}(x^{n+1})$
denotes a symmetric tensor on $\partial_{\infty}M$ whose components
in the $\left \{ \frac{\partial}{\partial y^A} \right \}$ basis are
each bounded above in magnitude by $Cx^{n+1}$ for some constant $C$
as $x\to 0$. We define the {\it mass} of such a metric to be
\begin{equation}
m\equiv m[g]:=\int_{\partial_{\infty}M}g_{(k)}^{AB}\kappa_{AB}d\mu_{(k)}
\ , \label{eq2.6}
\end{equation}
where $d\mu_{(k)}$ is the $g_{(k)}$ metric volume element on
$\partial_{\infty}M$, and
$g_{(k)}^{AB}\kappa_{AB}:\partial_{\infty}M\to {\mathbb R}$ is
called the {\it mass aspect function}.
\bigskip

\section{Term-by-term flow}
\setcounter{equation}{0}

\noindent Bahuaud \cite{Bahuaud} proves the short time existence of
a solution to the Ricci-DeTurck flow of asymptotically hyperbolic,
conformally compactifiable metrics $g(t)$ having the form
\begin{equation}
g(t)=\frac{dx^2+{\hat h}(x,y^A)+v(t,x,y^A)}{x^2}\ , \label{eq3.1}
\end{equation}
where $y^A$ are coordinates on $\partial_{\infty}M$ and ${\hat h}$
is a fixed, time independent metric on level sets of $x$, so that
$\frac{dx^2+{\hat h}}{x^2}=:h$ is the initial asymptotically
hyperbolic metric and $v(0,x,y^A)=0$.  This means that $x$ is a {\it
special defining function} for $h$, and $x$ is fixed in time during
Bahuaud's Ricci-DeTurck flow. The DeTurck diffeomorphism $\varphi_t$
is generated by a vector field $X=X^k\frac{\partial}{\partial x^k}=\frac{\partial \varphi_t}{\partial
t}$ which vanishes on the conformal boundary (see \cite{Bahuaud}
section 2, or appendix B for the explicit expression for $X$).
This means that the components of $X$ in the coordinate basis of
(\ref{eq3.1}) are ${\cal O}(x^2)$. Then the pullback metric is
\begin{eqnarray}
\varphi_t^* g(t)&=&\frac{\varphi_t^*\left ( dx^2+{\hat h} \right )
+\varphi_t^*v}{\left ( x\circ \varphi_t \right )^2}\nonumber\\
&=& \frac{dx^2+{\hat h}+v(t,x,y^A)+\left ( \varphi_t^*-{\rm id} \right )
\left ( dx^2+{\hat h}+v\right ) }{\left ( x\circ \varphi_t \right )^2} \nonumber\\
&=& \frac{dx^2+{\hat h}+ {\bar v}(t,x,y^A)}{\left ( x\circ \varphi_t \right )^2}
\ . \label{eq3.2}
\end{eqnarray}
where ${\bar v}:=v+\left ( \varphi_t^*-{\rm id} \right ) \left (
dx^2+{\hat h} +v \right )$. Since $X$ vanishes at $x=0$, then
$X^k=\frac{d}{dt}\left ( x^k\circ\varphi_t\right )\in{\cal O}(x^2)$
in our $(x^k)=(x,y^A)$ coordinate basis. Integrating this over $t\in
[0,T]$ with $\varphi_0={\rm id}$ yields
\begin{equation}
x^k\circ\varphi_t=x^k+{\cal O}(x^2)\ , \label{eq3.3}
\end{equation}
where $z\in{\cal O}(x^2)$ means that $|z|< c_Tx^2$ for $t\in [0,T]$,
where $c_T\in{\mathbb R}_+$ can depend on $T$. Since $x^1=x$, this
yields $x\circ \varphi_t = x+{\cal O}(x^2)$. Moreover,
differentiating (\ref{eq3.3}), then $\frac{\partial}{\partial x^j}
\left ( x^k\circ\varphi_t \right ) =\delta^k_j+{\cal O}(x)$, from
which it follows that $\left [\varphi_t^*-{\rm id}\right
]_j^k\in{\cal O}(x)$. Hence, since in addition $v\in{\cal O}(x)$,
then ${\bar v}\in{\cal O}(x)$.

We may therefore write the pullback metric which evolves under
(\ref{eq1.1}) as
\begin{eqnarray}
g(t)&=&\frac{1}{\rho_{(k)}^2(x)}{\tilde g}(t)
=\frac{1}{\rho_{(k)}^2(x)}\left [ dx^2+g_{(k)} +x^m\kappa(t)/m
+{\cal O}(x^{m+1})\right ] \ ,\nonumber \\
\rho_{(k)}(x)&=&\begin{cases}{\rm sinh}(x)\ ,&{\rm if\ }k=1\ ,\\
x\ ,&{\rm if\ }k=0\ ,
\end{cases}\label{eq3.4}\\
{\rm Ric}[g_{(k)}]&=&(n-2)kg_{(k)}\ , \nonumber
\end{eqnarray}
where we take $1\le m\le n =\dim M$ and where
$\kappa(t):=\kappa_{ij}(t)dx^idx^j$. This is more general than
(\ref{eq2.5}) as we do not assume that $m=n$ (indeed, we start our
iteration below with $m=1$) and we allow that that $\kappa_{11}$ and
$\kappa_{1A}$ can be nonzero ($\rho$ is not assumed to be a special
defining function at arbitrary $t$). Note that $\kappa_{ij}$ is a
function of the coordinates $y^A$ on level sets of $x$. The metric
induced from ${\tilde g}$ on $x=const$ hypersurfaces is
\begin{equation}
{\hat g}_{AB}:=g_{(k)AB}+x^m\kappa_{AB}/m+{\cal O}(x^{m+1})\ ,
\label{eq3.5}
\end{equation}
and we also define
\begin{equation}
{\tilde K}_{AB}:=\frac12 \frac{\partial}{\partial x}{\hat g}_{AB}
= \frac12 x^{m-1}\kappa_{AB}+{\cal O}(x^m)\ . \label{eq3.6}
\end{equation}
This is nearly the extrinsic curvature (but not quite, since
$\frac{\partial}{\partial x}$ need not be a unit vector). Although
the components of ${\hat g}_{AB}$ are simply ${\tilde g}\left (
\frac{\partial}{\partial y^A}, \frac{\partial}{\partial y^B}\right
)$, we find it somewhat useful to maintain a notational distinction
between ${\hat g}_{AB}$ and ${\tilde g}_{AB}$.

We wish to prove that if the initial metric has this form with
$m=n$, then so does the flowing metric at any time $t>0$ along the
flow. To do so, we must first compute the right-hand side of
(\ref{eq1.1}) for a metric of the above form.

We employ several standard expressions. The first is the radial {\it
matrix Riccati equation}:
\begin{equation}
\frac{\partial}{\partial x}{\tilde K}_{AB}
-{\tilde K}_{AC}{\tilde K}^C{}_B
= -{\tilde R}^1{}_{A1B}+{\cal O}(x^{2m-2})\ , \label{eq3.7}
\end{equation}
where ${\tilde R}^1{}_{A1B}=\left \langle dx, {\widetilde {\rm
Riem}}\left ( \frac{\partial}{\partial x}, \frac{\partial}{\partial
y^B} \right ) \frac{\partial}{\partial y^A}\right \rangle$ (the
``1'' indicates the $x$-direction), and the tilde denotes the
curvature tensor of ${\tilde g}=g/\rho_{(k)}^2$. The ${\cal
O}(x^{2m-2})$ correction term occurs only because $K_{AB}$
approximates the extrinsic curvature of $x=const$ surfaces only to
this order. The next ingredients are the equations of Gauss,
Codazzi, and Mainardi. If we write ${\tilde K}:={\tilde K}^A{}_A$
and if ${\hat R}_{AB}$ denotes the Ricci curvature of the connection
${\hat \nabla}$ of the induced metric ${\hat g}$, these yield
\begin{eqnarray}
{\tilde R}_{11}&=&{\tilde R}-{\hat R}+{\tilde K}^2
-{\tilde K}_{AB}{\tilde K}^{AB}\ ,
\label{eq3.8} \\
{\tilde R}_{1A}&=&{\tilde R}_{A1}={\hat \nabla}^B{\tilde K}_{AB}
-{\hat \nabla}_A {\tilde K} \label{eq3.9}\ , \\
{\tilde R}_{AB}&=& {\hat R}_{AB}+{\tilde R}^1{}_{A1B}
+{\tilde K}_{AC}{\tilde K}^C{}_B-{\tilde K}_{AB}{\tilde K}
\label{eq3.10}\ .
\end{eqnarray}
We also need that Taylor's theorem gives
\begin{equation}
{\hat R}_{AB}=R_{AB}[g_{(k)}]+{\cal O}(x^m)=(n-2)kg_{(k)AB}
+{\cal O}(x^m)\ . \label{eq3.11}
\end{equation}

The trace of (\ref{eq3.7}), together with (\ref{eq3.6}), allows us
to estimate ${\tilde R}_{11}$. We can estimate ${\tilde R}_{1A}$
immediately from equation (\ref{eq3.9}). To estimate ${\tilde
R}_{AB}$ we combine (\ref{eq3.7}) and (\ref{eq3.10}) and use
(\ref{eq3.6}). We obtain
\begin{eqnarray}
{\tilde R}_{11}&=&-\frac12 (m-1) x^{m-2}g_{(k)}^{AB}\kappa_{AB}
+{\cal O}(x^{m-1})\ , \label{eq3.12}\\
{\tilde R}_{A1}&=&{\tilde R}_{1A}={\cal O}(x^{m-1})\ , \label{eq3.13}\\
{\tilde R}_{AB}&=&(n-2)kg_{(k)AB}-\frac12 (m-1)x^{m-2}\kappa_{AB}
+{\cal O}(x^{m-1})\ , \label{eq3.14}
\end{eqnarray}
where, when $m=1$, these expressions are valid if we
interpret $(m-1)x^{m-2}$ as being identically zero, including at
$x=0$.

Another standard expression relates the Ricci curvature $R_{ij}$ of
$g$ to the Ricci curvature ${\tilde R}_{ij}$ of ${\tilde g}$. It is
\begin{equation}
R_{ij}={\tilde R}_{ij}+\frac{1}{\rho_{(k)}} \left [
(n-2) {\tilde \nabla}_i {\tilde \nabla}_j \rho_{(k)}
+{\tilde g}_{ij}{\tilde \Delta}\rho_{(k)} \right ]
-(n-1){\tilde g}_{ij}\frac{\widetilde {\vert \nabla \rho_{(k)}
\vert}^2}{\rho_{(k)}^2}\ , \label{eq3.15}
\end{equation}
where ${\tilde \Delta}:={\tilde g}^{ij}{\tilde \nabla}_i{\tilde
\nabla}_j$ and ${\widetilde {\vert \nabla \rho_{(k)}
\vert}^2}:={\tilde g}^{ij}{\tilde \nabla}_i \rho_{(k)}{\tilde
\nabla}_j \rho_{(k)}$. After straightforward computation using
expressions (\ref{eq3.12}--\ref{eq3.15}), we obtain
\begin{eqnarray}
R_{11}&=&-\frac{(n-1)}{\rho_{(k)}^2}-\frac12(n-1)x^{m-2}\kappa_{11}
-\frac12 (m-2)x^{m-2}g_{(k)}^{CD}\kappa_{CD}\nonumber\\
&&+{\cal O}(x^{m-1})\ , \label{eq3.16}\\
R_{1A}&=&R_{A1}=-\frac{(n-1)}{m}x^{m-2}\kappa_{1A}+{\cal O}(x^{m-1})
\ , \label{eq3.17}\\
R_{AB}&=&-\frac{(n-1)}{\rho_{(k)}^2}\left (kg_{(k)AB}
+\frac{1}{m}x^{m}\kappa_{AB} \right )\nonumber\\
&&+\frac12\left [ g_{(k)}^{CD}\kappa_{CD}+\left ( \frac{2n-2}{m}-1
\right ) \kappa_{11}\right ]x^{m-2}g_{(k)AB}\nonumber\\
&&+\frac12 \left ( n-m-1\right ) x^{m-2}\kappa_{AB}+{\cal O}
(x^{m-1})\ . \label{eq3.18}
\end{eqnarray}
Using $E_{ij}:=R_{ij}+(n-1)g_{ij}$, then we obtain
\begin{eqnarray}
E_{11}&=&
-\frac12 (m-2)\left [ \frac{(n-1)}{m}\kappa_{11}
+ g_{(k)}^{CD}\kappa_{CD}\right ]x^{m-2}+{\cal O}(x^{m-1})
\ , \label{eq3.19}\\
E_{1A}&=&E_{A1}={\cal O}(x^{m-1})\ , \label{eq3.20}\\
E_{AB}&=&\frac12\left [ \left ( \frac{2n-2}{m}-1\right )\kappa_{11}
+ g_{(k)}^{CD}\kappa_{CD}\right ]x^{m-2}g_{(k)AB}\nonumber\\
&&+\frac12 (n-m-1) x^{m-2} \kappa_{AB}+{\cal O}(x^{m-1})
\ . \label{eq3.21}
\end{eqnarray}

\bigskip
\noindent{\bf Proposition 3.1.} {\sl Let $g(t)$ be a solution of
(\ref{eq1.1}) of the form (\ref{eq3.4}) on some interval
$t\in[0,T)$. If the expansion (\ref{eq3.4}) for the initial metric
$g(0)$ begins at order $m=n$, then so does the expansion for
$g(t)$.}

\bigskip
\noindent{\bf Proof.} From (\ref{eq3.4}), $\frac{\partial
g_{ij}}{\partial t} = \frac{1}{m}x^{m-2}\frac{\partial
\kappa_{ij}}{\partial t}+{\cal O}(x^{m-1})$. Then (\ref{eq1.1})
yields
\begin{equation}
\frac{\partial \kappa_{ij}}{\partial t}=-\frac{2m}{x^{m-2}}
E_{ij}+{\cal O}(x)
= A_{ij}{}^{kl}\kappa_{kl}+{\cal O}(x)\ , \label{eq3.22}
\end{equation}
where $A$ is a matrix with components given (using
(\ref{eq3.19}--\ref{eq3.21})) by
\begin{eqnarray}
A_{11}{}^{11}&=&(m-2)(n-1) \ , \label{eq3.23}\\
A_{11}{}^{CD}&=&m (m-2)g^{CD}_{(k)} \ , \label{eq3.24}\\
A_{AB}{}^{11}&=&-m\left ( \frac{2n-2}{m}-1 \right )g_{(k)AB}
\ , \label{eq3.25}\\
A_{AB}{}^{CD}&=&-m g_{(k)AB}g_{(k)}^{CD}
-m (n-m-1)\delta_A^C\delta_B^D\ . \label{eq3.26}
\end{eqnarray}
Taking the limit $x\to 0$, we obtain a linear system for
$\kappa_{ij}(t)$. If $m<n$, then $\kappa_{ij}(0)=0$ by assumption.
But then, by uniqueness, $\kappa_{ij}(t)=0$. \qed

\bigskip

\section{Ricci flow and mass}
\setcounter{equation}{0}

\noindent We now wish to consider the Ricci flow developing from an
initial metric of the form (\ref{eq2.5}). This means that $m=n$ and
that $\kappa_{11}(0)=\kappa_{A1}(0)=0$. Proposition 3.1 then
guarantees that $m=n$ throughout the flow, but it does not guarantee
that $\kappa_{11}(t)=0$, nor that $\kappa_{A1}(t)=0$, at any $t>0$.
Since we wish to compute the mass of the flowing metric, we would
have to transform this metric back to the form (\ref{eq2.5})
(vanishing $\kappa_{i1}$) whenever we wish to compute the mass.
Alternatively, we can transform the mass formula to a form valid for
arbitrary $\kappa_{ij}$.

To this end, consider a metric of the form
\begin{eqnarray}
g&=&\frac{1}{\rho_{(k)}^2(x')}\left [ dx'^2+g_{(k)}
+\frac{x'^n}{n}\kappa_{ij}dx^idx^j+{\cal O}(x'^{n+1})\right ]\nonumber\\
&=&\frac{1}{\rho_{(k)}^2(x')}\bigg [ dx'^2+g_{(k)} +\frac{x'^n}{n}\left
( \kappa_{11}dx'^2+2\kappa_{1A}dx'dz^A+\kappa_{AB}dz^Adz^B \right )
\nonumber\\
&&+{\cal O}(x'^{n+1})\bigg ] \ , \label{eq4.1}
\end{eqnarray}
generalizing the form (\ref{eq2.5}) to allow for possibly nonzero
$\kappa_{11}$ and $\kappa_{1A}$. The coordinate transformation
\begin{equation}
y^A=z^A+\frac{x'^{n+1}}{n(n+1)}g_{(0)}^{AB}\kappa_{1B}\ . \label{eq4.2}
\end{equation}
brings (\ref{eq4.1}) to the form
\begin{equation}
g=\frac{1}{\rho_{(k)}^2(x')}\left [ dx'^2+g_{(k)} +\frac{x'^n}{n}\left
( \kappa_{11}dx'^2+\kappa_{AB}dy^Ady^B \right )
+{\cal O}(x'^{n+1})\right ] \ , \label{eq4.3}
\end{equation}
The further transformation
\begin{equation}
x=x'\left ( 1 + \frac{1}{2n^2}\kappa_{11}x'^n\right )\label{eq4.4}
\end{equation}
brings the metric to the form
\begin{equation}
g=\frac{1}{\rho_{(k)}^2(x)}\left [ dx^2+g_{(k)} +\frac{x^n}{n}\left
( \kappa_{AB}+\frac{1}{n}\kappa_{11}g_{(k)AB} \right )dy^Ady^B
+{\cal O}(x^{n+1})\right ] \ . \label{eq4.5}
\end{equation}
Once we have this form, Lemma 3.10 of \cite{ACG} shows that further
transformations eliminate the higher order terms in $g_{11}$ and
$g_{1A}$, returning $g$ to the form (\ref{eq2.5}) but with
$\kappa_{AB}$ replaced by
$\kappa_{AB}+\frac{1}{n}\kappa_{11}g_{(k)AB}$. Thus we have proved
the following:

\bigskip
\noindent{\bf Lemma 4.1.} {\sl The mass of a metric of the form
(\ref{eq4.1}) is}
\begin{equation}
m\equiv m[g]:=\int_{\partial_{\infty}M}\left (g_{(k)}^{AB}
\kappa_{AB}+\frac{(n-1)}{n}\kappa_{11}\right ) d\mu_k
\ . \label{eq4.6}
\end{equation}
\bigskip

\noindent We will give the integrand in the above formula a name.

\bigskip
\noindent{\bf Definition 4.2.} The {\it mass aspect function}
$\sigma:\partial_{\infty}M\to {\mathbb R}$ is
\begin{equation}
\sigma:=g_{(k)}^{AB} \kappa_{AB}
+\frac{(n-1)}{n}\kappa_{11} \ . \label{eq4.7}
\end{equation}

\noindent{\bf Proposition 4.3.} {\sl Let $g(t)$ be as in proposition
3.1, with the expansion (\ref{eq3.4}) for $g(0)$ beginning at order
$m=n$. Then the mass aspect of $g(t)$ evolves as
\begin{equation}
\sigma(t)=\sigma_0 e^{-(n-2)t}\ , \label{eq4.8}
\end{equation}
where $\sigma_0=\sigma(0)$.}

\bigskip
\noindent{\bf Proof.} By proposition 3.1, the metric takes the form
(\ref{eq4.1}) at each $t\in [0,T)$. The flow of $\kappa_{ij}(t)$ is
then obtained by setting $m=n$ in (\ref{eq3.22}--\ref{eq3.26}). This
gives
\begin{eqnarray}
\frac{\partial \kappa_{11}}{\partial t} &=& (n-2) \left [ (n-1)
\kappa_{11} +ng^{CD}_{(k)}\kappa_{CD} \right ] \ , \label{eq4.9}\\
\frac{\partial \kappa_{AB}}{\partial t} &=& -(n-2)\kappa_{11}g_{(k)AB}
+n\left [ \kappa_{AB} -g_{(k)AB}g^{CD}_{(k)}\kappa_{CD}\right ]
\ . \label{eq4.10}
\end{eqnarray}
Forming the appropriate linear combination, we thus obtain
\begin{equation}
\frac{\partial \sigma}{\partial t} = -(n-2)\sigma\ , \label{eq4.11}
\end{equation}
from which the proposition follows. \qed
\bigskip

Theorem 1.1 now immediately follows from Proposition 4.3.

\bigskip
\noindent{\bf Proof of Theorem 1.1.} Integrate (\ref{eq4.8}) over
$\partial_{\infty}M$.\qed
\bigskip

In passing, we also obtain a theorem known from \cite{AnderssonDahl}
for $k=1$ (with $\partial_{\infty}M$) and generalize it to the $k=0$
case.

\bigskip
\noindent{\bf Lemma 4.4.} {\sl If an asymptotically hyperbolic
metric is Einstein, it cannot have nonzero mass.}

\bigskip
\noindent{\bf Proof.} Setting $m=n$ in (\ref{eq3.19}) yields
$E_{11}=-\frac12 (n-2)\sigma$. The Einstein condition is $E_{ij}=0$,
so in particular $E_{11}=0$. Thus, the mass aspect vanishes. \qed
\bigskip

\noindent When $\partial_{\infty} M\cong S^n$, a partial converse of
this result is proved in the next section.

\section{Rigidity}
\setcounter{equation}{0}

\noindent In this section, we give a variant of a standard proof
\cite{ACG} of Proposition 1.4 based on the flow (\ref{eq1.1}) and
following ideas in \cite{Haslhofer}. The variant method relies on
Bahuaud's short-time existence \cite{Bahuaud} theorem to construct a
metric variation. By way of contrast, \cite{ACG} directly posits a
suitable variation (the first Newton approximation to a solution of
(\ref{eq1.1})). In return for this added complexity at one point in
the argument, our method gains in simplicity at another by avoiding
the necessity in \cite{ACG} (section 3.2.3 thereof) of replacing a
solution of the Yamabe equation by a solution of its linearization.

To begin, we say that the mass aspect of an asymptotically
hyperbolic metric $g$ is {\it of semi-definite sign} if
$\sigma:\partial_{\infty}\to {\mathbb R}$ is either a nonnegative
function or a nonpositive function; i.e., if $\sigma(p)\sigma(q)\ge
0$ $\forall p,q\in \partial_{\infty}M$. We recall from Proposition
1.4 the class ${\cal G}$ of metrics whose elements are
asymptotically hyperbolic with $\partial_{\infty}M\simeq S^{n-1}$
and with well-defined mass, obeying $E[g]:=R+n(n-1)\ge 0$, $3\le
n=\dim M \le 6$, and further restricted if $M$ is not spin to
contain only metrics whose mass aspect is of semi-definite sign. We
take $3\le n=\dim M \le 6$.

\bigskip
\noindent {\bf Theorem 5.1 (Positive Energy Theorem \cite{Wang, CH,
ACG}).} {\sl If $g\in{\cal G}$ then $m[g]\ge 0$.}

\bigskip
\noindent {\bf Proof.} For the spin case, see \cite{Wang, CH}. For
the semi-definite mass aspect case, see \cite{ACG}. \qed

\bigskip
\noindent {\bf Proof of Proposition 1.4.} By way of contradiction,
assume that $g\in{\cal G}$ and $E_{ij}\vert_p\neq 0$ at some $p\in
M$. Without loss of generality, we can take $E\vert_p>0$, for if it
is not, then we can take $g$ to be initial data for a flow $g(t)$
solving (\ref{eq1.1}). By \cite{Bahuaud} and Proposition 3.1, $g(t)$
will exist at least on a non-empty interval $[0,T_1)$ and will
remain asymptotically hyperbolic with well-defined mass. Then a
standard derivation shows that under the flow (\ref{eq1.1}),
$E:=R+n(n-1)$ evolves according to
\begin{equation}
\frac{\partial E}{\partial t} =\Delta E +2 E_{ij}E^{ij}-2(n-1)E
\ . \label{eq5.1}
\end{equation}
By the maximum principle and since we assume that $E(0)\ge 0$, then
$E(t)\ge 0$ $\forall t\in [0,T_1)$. Furthermore, inspection of
(\ref{eq5.1}) shows that if $E(0)=0$ and $E_{ij}(0)\vert_p\neq 0$,
then $\frac{\partial E}{\partial t}(0)\vert_p>0$, and so
$E(t)\vert_p>0$ for $t\in (0,T)$ for some $0<T<T_1$. Furthermore, by
Theorem 1.1, $m(t)=0$.

Hence, we now have $E\ge 0$ on $M$ and a $p\in M$ such that
$E\vert_p>0$, and we have $m=0$. When $k=1$, Proposition 3.13 of
\cite{ACG} now shows that there is a conformal transformation of $g$
which produces a new metric ${\hat g} = w^{\frac{4}{n-2}}g$, $w\to
1$ on approach to $\partial_{\infty}M$, such that ${\hat g}$ is
strongly asymptotically hyperbolic, the scalar curvature of ${\hat
g}$ obeys ${\hat E}:={\hat R}+n(n-1)=0$, and the mass aspect of
${\hat g}$ is pointwise strictly lower than that of $g$: ${\hat
\sigma}<\sigma$. The proof of this proposition relies on the
existence of a suitable solution of the Yamabe prescribed scalar
curvature equation on $(M,g)$, which was proved in \cite{AM}. Of
course, it then follows that the mass of ${\hat g}$ is strictly
lower than that of $g$: ${\hat m}:=m({\hat g})<m=m(g)=0$. But if
$g\in{\cal G}$, then ${\hat g}\in{\cal G}$, since ${\hat E}=0$ and
since the flow and conformal transformation preserve strong
asymptotic hyperbolicity. Furthermore, if $(M,g)$ is spin then so is
$(M,{\hat g})$ since possession of a spin structure is a topological
property, and if $\sigma$ is of semi-definite sign, then
$m=0\Rightarrow \sigma=0$, and so ${\hat \sigma}<\sigma$ implies
that ${\hat \sigma}<0$. Either way, this contradicts the positive
mass theorem (\cite{Wang, CH, ACG}). Hence $E_{ij}\equiv 0$
pointwise on $M$. Then, using the Einstein manifold rigidity result
of \cite{Qing} valid for $3\le n=\dim M \le 6$, $(M,g)$, we see that
$(M,g)$ is isometric to standard hyperbolic space. \qed

\bigskip

\section{Parabolic scaling and ADM mass}
\setcounter{equation}{0}

\noindent The Ricci flow of mass has arisen in the physics of string
theory and the renormalization group. In that context, \cite{GHMS}
exhibits an explicit example of a 2-dimensional Ricci soliton whose
ADM mass (in 2-dimensions, this is the deficit angle of the cone to
which the metric asymptotes) remains constant during Ricci flow. The
flow has a $t\to\infty$ limit, which is flat space. Thus the mass
remains constant during the flow and equal to the mass of the
initial metric, but the mass of the limit manifold is zero. This
behaviour persists in much more general circumstances. In \cite{DM,
OW} it was shown that ADM mass is constant under asymptotically flat
Ricci flow in all dimensions $n\ge 3$, and in \cite{OW} it was shown
that rotationally symmetric, asymptotically flat Ricci flow
approaches flat space in the $t\to\infty$ limit if no minimal
surface is present in the initial data.

A short heuristic argument which we
now outline shows that the behaviour of mass in the asymptotically
flat flow can be ``derived'', based on  reasoning that follows from
(\ref{eq1.3}). When the asymptotic sectional curvature is $1/\ell^2$
and is no longer normalized to $1$, parabolic scaling leads us to
replace $t\mapsto t/\ell^2$ in (\ref{eq1.3}):
\begin{equation}
m(t):=m_0 e^{-(n-2)t/\ell^2}\ . \label{eq1.4}
\end{equation}
This is the formula for the flow of mass under the evolution
\begin{equation}
\frac{\partial g_{ij}}{\partial t} = -2\left ( R_{ij}
+\frac{(n-1)}{\ell^2} g_{ij}\right )\ , \label{eq1.5}
\end{equation}
which is obtained from (\ref{eq1.1}) by the replacements $t\mapsto
t/\ell^2$ and $g\mapsto g/\ell^2$. If we now consider the limit
$\ell\to\infty$ at any fixed $t$, we see that (\ref{eq1.4}) becomes
simply $m(t)=m_0=const$, (\ref{eq1.5}) becomes the usual Ricci flow,
and since the asymptotic radius of curvature $\ell$ goes to infinity
the solution $g$ becomes asymptotically flat. However, if instead we
take $t\to\infty$ first, before sending $\ell\to\infty$, the mass would
still of course approach $0$ in this limit. While this discussion is
quite far from rigourous, we see that it reproduces the known
behaviour of the mass of an asymptotically flat metric under Ricci
flow.

Finally, we note that \cite{GHMS} conjectured that an initially
positive mass will decrease under a string theory process called
{\it tachyon condensation}, which can be modelled by Ricci flow.
Asymptotically flat spacetimes can have two notions of mass, called
ADM and Bondi \cite{Wald}. The conjecture referred directly to the
Bondi mass, but the Ricci flow occurs on a Riemannian manifold where
there is only ADM mass, whose flow behaviour is at best a very
coarse manifestation of the conjectured behaviour. When one passes
to asymptotically anti-de Sitter spacetimes, the Lorentzian analogue
of asymptotically hyperbolic manifolds, the distinction between
these two notions of mass disappears. Perhaps this is why the smooth
mass evolution (\ref{eq1.3}) under (\ref{eq1.1}) seems to reflect so
well the expectation arising from physics of tachyon condensation.

\appendix
\section{Equivalence of Wang and Chru\'sciel-Herzlich masses}
\setcounter{equation}{0}

\noindent In his paper \cite{Wang}, Wang defines the mass
(\ref{eq2.6}) for $k=1$ and $\partial_{\infty}M\simeq S^{n-1}$. This
mass is equivalent to the Chru\'sciel-Herzlich mass \cite{CH}, which
is formulated in a chart-independent manner. Let $b$ and $g$ be any
two Riemannian metrics on a manifold $M^n$. Let $D$ be the
Levi-Civit\`a connection of $b$ and let $\nabla$ be the
Levi-Civit\`a connection of $g$. The \emph{Chru\'{s}ciel-Herzlich
mass} of $(M^n,g)$ with respect to the \emph{reference} metric $b$
is
\begin{eqnarray}
m_{\rm CH}(g) &=& \lim_{R\to\infty}\int_{N_R}\mathbb{U}^i
(V\circ \phi^{-1})dS_i \ , \label{eqA.1}\\
\mathbb{U}^i &=& \sqrt{g}\left\{V\left (
g^{ik}g^{jl}-g^{ij}g^{kl}\right ) )D_jg_{kl}
+\left (g^{jk}\nabla^iV-g^{ik}\nabla^jV\right )
e_{jk}\right\}\ ,\label{eqA.2}\\
e_{jk}&:=& g_{jk}-b_{jk}\ , \label{eqA.3}
\end{eqnarray}
where $V:M_{\rm ext}\to\mathbb{R}$ is a smooth ${\cal O}(r)$
function on a collar neighbourhood $M_{\rm ext}$ of
$\partial_{\infty}M$, $\phi^{-1}:M_{ext}\to [R,\infty)\times N$ is a
smooth diffeomorphism, and $dS_i := n_i\,dA_R$ where $n_i$ is the
unit normal and $dA_R$ is the induced volume (i.e., area) element on
$N_R$.

In (\ref{eq2.6}) we extended Wang's definition to the case of $k=0$,
with $\partial_{\infty}M$ a flat torus. The Chru\'sciel-Herzlich
mass \cite{CH} formulation already includes this case. Here we show
that our $k=0$ extension agrees with the chart-independent $k=0$
Chru\'sciel-Herzlich formulation, as required for the proof of
Corollary 1.3. We will do this by evaluating the
Chru\'sciel-Herzlich mass in the coordinate gauge in which the
metric takes the form (\ref{eq2.5}) (so that
$\kappa_{1i}=\kappa_{i1}=0$) and checking that the
Chru\'sciel-Herzlich expression for the mass reduces to
(\ref{eq2.6}) with $k=0$. The $k=1$ case is similar and
straightforward, and has been reported in \cite{CH}.

In the $k=0$ case, for which $N$ is an $(n-1)$-torus, we have (up to
a scale for $V$, which we fix to be $1$)
\begin{eqnarray}
V(r)&=&r\ , \label{eqA.4}\\
b &=& \frac{dr^2}{r^2}+r^2g_{(0)}=\frac{dr^2}{r^2}
+r^2\sum_{A=2}^n d\theta_A^2\ , \label{eqA.5}\\
e_{jk} &=& g_{jk}-b_{jk} = \frac{1}{nr^{n-2}}\kappa_{jk}
+{\cal O}(1/r^{n-1})\ , \label{eqA.6}
\end{eqnarray}
where in the last equation we used the form (\ref{eq2.5}) for $g$
with $k=0$ and $\rho_{(0)}(r)=r=1/x$, and in the sequel we will use
$A,B\in \{2,\dots,n\}$ to denote coordinates on hypersurfaces
$r=const.$

It is now a matter of plugging (\ref{eq2.5},
\ref{eqA.4}--\ref{eqA.6}) into (\ref{eqA.2}). We write
\begin{eqnarray}
\mathbb{U}^i &=& \sqrt{g}\, r \left\{\mathbb{A}^i+\mathbb{B}^i
\right\} \ , \label{eqA.7}\\
\mathbb{A}^i &=& \left(g^{ik}g^{jl}-g^{ij}g^{kl}\right)D_jg_{kl}
\ , \label{eqA.8}\\
\mathbb{B}^i &=& \left(\frac{\nabla^i r}{r}g^{jk}
-\frac{\nabla^j r}{r}g^{ik}\right)e_{jk} \nonumber\\
&=&\left(g^{il}g^{jk}-g^{jl}g^{ik}\right)\frac{\nabla _l r}{r}
e_{jk}\ ,\label{eqA.9}
\end{eqnarray}
and find that
\begin{eqnarray}
{\mathbb A}^1&=&\frac{(n-1)}{nr^{n-1}}g_{(0)}^{AB}\kappa_{AB}
+{\cal O}(1/r^n)\ , \label{eqA.10}\\
{\mathbb B}^1&=&\frac{1}{nr^{n-1}}g_{(0)}^{AB}\kappa_{AB}
+{\cal O}(1/r^n)\ , \label{eqA.11}
\end{eqnarray}
and ${\mathbb A}^A={\cal O}(1/r^n)$, ${\mathbb B}^A={\cal
O}(1/r^n)$. Substituting (\ref{eqA.10}) and (\ref{eqA.11}) into
(\ref{eqA.7}) yields
\begin{equation}
{\mathbb U}^i = \sqrt{g} \left\{\frac{1}{r^{n-2}}g_{(0)}^{AB}
\kappa_{AB}\delta^i_1+{\cal O} (1/r^{n-1})\right\} \ . \label{eqA.12}
\end{equation}
Plugging this into (\ref{eqA.1}) yields
\begin{equation}
m_{\rm CH}(g) = \lim_{R\to\infty}\int_{N_R}\mathbb{U}^i\,dS_i
=\int_{N}g_{(0)}^{AB}\kappa_{AB}\,dV(g_{(0)})\ , \label{eqA.14}
\end{equation}
which is the mass formula (\ref{eq2.6}) for the $k=0$ case. As
discussed in the Introduction, for the metric (\ref{eq1.2}) this
formula yields a mass of $-\frac{4\pi}{n}\prod_{i=3}^na_i$.

\section{Ricci-DeTurck flow}
\setcounter{equation}{0}

\noindent
Bahuaud's normalized Ricci-DeTurck flow \cite{Bahuaud} is
\begin{equation}
\frac{\partial g}{\partial t}=-2E+\pounds_Xg
\equiv -2\left ( {\rm Ric}+(n-1)g \right )+\pounds_Xg\ , \label{eqB.1}
\end{equation}
where $\pounds_Xg$ is the Lie derivative of $g$ along the DeTurck vector field $X$. Pulling back along the diffeomorphsim generated by $X$ yields the normalized Ricci flow (\ref{eq1.1}).
One can apply the iteration of section 3 directly to (\ref{eqB.1}), rather than to (\ref{eq1.1}) as is done in the main text.
We outline the
main steps. For comparison purposes, we fix $k=1$ so that $\partial
M$ carries a metric $g_{(1)}=g(S^{n-1},{\rm can})$ with sectional
curvature $1$, and adopt Bahuaud's notation
\begin{equation}
g=h+v\ , \label{eqB.2}
\end{equation}
where $g(0)=h$ so that  $v(0)=0$, $\frac{\partial h}{\partial t}=0$,
and the components of $v$ in an orthonormal basis of $h$ are ${\cal
O}(x)$. Since $h$ is now the initial metric, we write (using $\sinh\rho=\frac{4x}{4-x^2}$)
\begin{equation}
h:=x^{-2}{\tilde h}:=x^{-2}\left [ dx^2 +\left ( 1-(x/2)^2\right )^2
g(S^{n-1},{\rm can})+\frac{1}{n}x^n\kappa_{AB}\right ] \ , \label{eqB.3}
\end{equation}
where $x$ is a special defining function for the asymptotically hyperbolic
metric $h$. Note that $\kappa$ here corresponds to the initial
metric ``perturbation'' and is constant along the flow, unlike in the
main text. We write
\begin{equation}
v(t)=\frac{1}{m}w_{ij}(t)x^{m-2}\ , \label{eqB.4}
\end{equation}
where $w(t,x)=w_0(t)+{\cal O}(x)$. We then consider the cases
$m\le n$, using that $w(0)=0$ when $m\le n$.

Let $\nabla$ be the Levi-Civit\`a connection of $g$ with connection
coefficients $\Gamma^k_{ij}$ and let ${\mathring \nabla}$ be the
Levi-Civit\`a connection of $h$ with connection coefficients
${\mathring \Gamma}^k_{ij}$.
The DeTurck vector field $X=X^k\frac{\partial}{\partial x^k}$ has components
\begin{equation}
X^k = g^{ij}\left ( \Gamma^k_{ij}-
{\mathring \Gamma}^k_{ij} \right )
=h^{kl}h^{ij} \left ( {\mathring \nabla}_i v_{jl}
-\frac12 {\mathring \nabla}_l v_{ij} \right )
+{\cal O}(v^2,v\cdot\partial v)\ . \label{eqB.5}
\end{equation}
A calculation then yields
\begin{eqnarray}
X_1&=&-\frac{(n-2)}{m}w_{11}x^{m-1}+\left ( \frac{m-2}{2m} \right ) \left ( w_{11}-g_{(1)}^{AB}w_{AB}\right )x^{m-1}+{\cal O}(x^m)
\ , \label{eqB.6}\\
X_A&=&\left ( \frac{m-n}{m}\right )w_{1A}x^{m-1}+{\cal O}(x^m)
\ . \label{eqB.7}
\end{eqnarray}

We then compute that
\begin{eqnarray}
\left (\pounds_X g\right )_{11}
&=& -2(n-2)w_{11}x^{m-2}+(m-2) \left ( w_{11}-g_{(1)}^{AB}w_{AB}\right )x^{m-2}\nonumber\\
&&+{\cal O}(x^{m-1})\ , \label{eqB.8}\\
\left (\pounds_X g\right )_{1A}&=&\left (\pounds_X g\right )_{1A}
=\frac{1}{m}(m-n)(m+1)w_{1A}x^{m-2}+{\cal O}(x^{m-1})
\ , \label{eqB.9}\\
\left (\pounds_X g\right )_{AB}&=& \left [ \frac{2(n-2)}{m}
w_{11}x^{m-2}-\frac{(m-2)}{m}\left ( w_{11}-g_{(1)}^{CD}w_{CD}\right )x^{m-2}\right ] g_{(1)AB}\nonumber\\
&&+{\cal O}(x^{m-1})\ . \label{eqB.10}
\end{eqnarray}

For $m<n$, the evolution equations for the metric perturbation under
DeTurck flow are obtained by replacing $\kappa$ by $w$ in equations
(\ref{eq3.19}--\ref{eq3.21}) and adding to these the appropriate Lie
derivative component from equations (\ref{eqB.8}--\ref{eqB.10}) to
form $\varepsilon_{ij}:=E_{ij}-\frac12 \pounds_X g_{ij}$. Then
(\ref{eq3.22}) is replaced by
\begin{equation}
\frac{\partial w_{ij}}{\partial t}
=-\frac{2m}{x^{m-2}}\varepsilon_{ij}+{\cal O}(x)=:{\cal A}_{ij}{}^{kl}w_{kl}
+{\cal O}(x)\ . \label{eqB.11}
\end{equation}
Then (\ref{eqB.11}) is a linear system with zero initial data (cf
equation (\ref{eq3.22}), and hence $w_{ij}(t)=0$ for $m<n$, so
Proposition 3.1 holds for Bahuaud's Ricci-DeTurck flow.

For $m=n$, we must replace $\kappa$ in (\ref{eq3.19}--\ref{eq3.21})
by $\kappa+w$ (but not in (\ref{eqB.8}--\ref{eqB.10})) and proceed as
above, adding (\ref{eqB.8}--\ref{eqB.10}) to $E$. In particular, we
are interested in the evolution of the combination
$\frac{(n-1)}{n}w_{11}+g_{(1)}^{AB}w_{AB}$, as this is the
time-evolving part of the mass aspect. But, setting $m=n$ in
(\ref{eqB.8}) and (\ref{eqB.10}), we obtain
\begin{equation}
\frac{(n-1)}{n}\left (\pounds_X g\right )_{11}
+g_{(1)}^{AB}\left (\pounds_X g\right )_{AB}={\cal O}(x^{n-1})
\ . \label{eqB.12}
\end{equation}
Thus, the Lie derivative term in (\ref{eqB.1}) does not contribute
to the evolution of the mass aspect. Then it is easy to see that
Proposition 4.3 applies also to Bahuaud's Ricci-DeTurck flow, and
then so does Theorem 1.1.

Finally, now we may put $m=n$ in (\ref{eqB.6}) and (\ref{eqB.7})
and raise the index,
obtaining $X^k\in{\cal O}(x^{n+1})$. Thus $\varphi^k_t={\rm
id}+{\cal O}(x^n)$, so the mass is invariant under the pullback of
$g$ by $\varphi^k_t$ (see, for example, Theorem 3.4 of
\cite{Herzlich}). This reproduces the Ricci flow result in the main
text.

\end{document}